
\magnification=\magstep1
\input amstex
\documentstyle{amsppt}
\leftheadtext{E. Makai, Jr.$^{*}$, H. Martini}
\rightheadtext{Centrally symmetric convex bodies}
\topmatter
\title {Centrally symmetric convex bodies and
sections having maximal quermassintegrals
{\centerline{\rm{Stud. Sci. Math. Hungar. 49 (2) (2012), 189-199}}}}
\endtitle
\author E. Makai, Jr.$^{*}$, H. Martini$^{**}$
\vskip.2cm
{\centerline{\rm{DOI: 10.1556/SScMath.49.2012.2.1197}}}
\endauthor
\address 
E. Makai, Jr.$^{*}$,
\newline
Alfr\'ed R\'enyi Mathematical Institute, 
Hungarian Academy of Sciences,
\newline
H-1364 Budapest, Pf. 127, 
HUNGARY,
\newline
{\rm{http://www.renyi.mta.hu/\~{}makai}}
\newline
H. Martini$^{**}$,
\newline
Technische Universit\"at Chemnitz,
Fakult\"at f\"ur Mathematik,
\newline
D-09107 Chemnitz, GERMANY,
\newline
www.tu-chemnitz.de/mathematik/geometrie
\endaddress
\email makai.endre\@renyi.mta.hu, martini\@mathematik.tu-chemnitz.de
\endemail
\thanks $^*$ Research (partially) supported by 
Hungarian National Foundation for 
Scientific Research, grant nos. K68398, K75016, K81146, and
DFG-GZ MA 1329/26-1.
\newline
$^{**}$ Research partially supported by DFG.
\endthanks
\subjclass {\it 2000 Mathematics Subject Classification.} 
Primary: 52A20. Secondary: 52A39\endsubjclass

\abstract 
Let $d \ge 2$, and let $K \subset {\Bbb{R}}^d$ be a convex body containing the
origin $0$ in its interior. In a previous paper we have
proved the following. The body $K$ is $0$-symmetric if and only if the
following holds. For each $\omega \in S^{d-1}$, we have that the
$(d-1)$-volume of the intersection of $K$ and an arbitrary hyperplane, with
normal $\omega$, attains its maximum if the hyperplane
contains $0$. An analogous theorem, for $1$-dimensional sections and
$1$-volumes, has been proved long ago by Hammer (\cite{H}). In this paper
we deal with the ($(d-2)$-dimensional) surface area, or with lower dimensional
quermassintegrals of these intersections, and prove an analogous,
but local theorem, for small $C^2$-perturbations, or $C^3$-perturbations of
the Euclidean
unit ball, respectively.
\endabstract
\endtopmatter\document

\head 1. Introduction\endhead

Let $d \ge 2$, and let $K \subset {\Bbb{R}}^d$ be a {\it{centered}}, i.e.,
$0$-symmetric convex body. We have observed in \cite{MM\'O}, Problem 3.10, that
by
the Aleksandroff-Fenchel inequalities (cf., e.g., \cite{S}) we have the
following statement. Let $0 \le l \le d-2$ be an integer, let $\omega \in
S^{d-1}$,
let $t \in {\Bbb{R}}$, 
and let $\omega^\bot $ be the orthocomplement of $\omega$ in
${\Bbb{R}}^d$. Then the quermassintegrals
$$
W_l \left[(K \cap (\omega^\bot + t \omega)) - t \omega\right]\,,
$$
considered in $\omega^\bot $, attain their maxima for $t = 0$. In the same
Problem 3.10, we have posed the question, whether the converse implication
holds.
For $l = 0$, i.e., for the case of $(d-1)$-volume, we proved this converse
implication, cf. \cite{MM\'O}, Corollary 3.2.

In this paper, we deal with the cases $1 \le l \le d-2$, and prove an
infinitesimal variant of the converse implication, for small
$C^2$-perturbations of the 
Euclidean unit ball
for $l=1$, and for small $C^3$-perturbations of the Euclidean unit ball for
$2 \le l \le d-2$.


\head 2. Preliminaries \endhead

We write ${\Bbb{R}}^d$ 
for the $d$-dimensional Euclidean space, and $S^{d-1}$ for
its unit sphere, where $d \ge 2$. The origin is denoted by $0$. We write $W_i$
for the
$(k-i)$-dimensional quermassintegrals of convex bodies in affine
$k$-subspaces of ${\Bbb{R}}^d$ (\cite{BF}, \cite{S}).

Basically we use the notations of \cite{MM\'O}. Variable points of $S^{d-1}$ are
denoted by $\omega, \xi, \eta $.
We use polar coordinates on $S^{d-1}$. That is,
for some
$\xi \in S^{d-1}$, that we consider as the north pole, and for $\omega \in
S^{d-1}$, we write
$$
\omega = \xi \sin \psi + \eta \cos \psi, {\text{ where }} \, \eta \in
\xi^{\bot }\cap
S^{d-1}, {\text{ and }} \, -\pi/2 \le \psi \le \pi/2\,.
$$
Thus, $\xi^\bot \cap S^{d-1}$ is the {\it{equator}} and $\psi$ is the
{\it{geographic latitude}}, that will be more convenient to us than the
customarily used $\varphi = \pi/2-\psi$.
Then we write
$$
\omega = (\eta, \psi)\,.
$$
In particular,
$$
(\eta,0) = \eta\,.
$$
A function $f: = S^{d-1} \to {\Bbb{R}}$ is {\it{even}}, or {\it{odd}}, 
if, for all
$\omega \in S^{d-1}$, we have $f(-\omega) = f(\omega)$, or $f(-\omega)
= - f(\omega)$, respectively.

In ${\Bbb{R}}^d$ 
we will use polar coordinates $\omega,\varrho$, with $\omega \in
S^{d-1}$, and $\varrho \in [0,\infty)$ (i.e., the point $\omega \varrho \in
{\Bbb{R}}^d$ has polar
coordinates $\omega,\varrho$). Also, for fixed $\xi \in S^{d-1}$, we will use
cylindrical coordinates $\eta,r,t$, with $\eta \in \xi^\bot \cap S^{d-1}$,
and $r
\in [0,\infty)$
(together polar coordinates in $\xi^\bot$), and $t \in {\Bbb{R}}$. 
Here, for $x \in
{\Bbb{R}}^d$, we have
$$
t = \langle x,\xi \rangle \,, {\text{ and }} \, x = r \eta + t \xi\,.
$$
For $x \in {\text{bd}}\,K$, we will also write, in cylindrical coordinates,
$$
x = r(\xi, \eta, t) \eta + t \xi\,,
$$
where the first variable of $r$ refers to $\xi$, and the last variable means
that we consider the radial function of the intersection
$K \cap (\xi^\bot + t \xi)$, with respect
to the ``origin'' $t \xi$.

We have, for $x \in {\Bbb{R}}^d$, that
$$
\varrho \cos \psi = r, \, {\text{ and }} \, \varrho \sin \psi = t\,.
\tag 1
$$
Differentiating these formulas with respect to $\psi$, and then setting $\psi
= 0$, we obtain
$$
\frac{\partial r}{\partial \psi}{\big|}_{\psi =0}
= \frac{\partial \varrho}{\partial
\psi}{\big|}_{\psi=0}, \, {\text{ and }} \,
\frac{\partial t}{\partial \psi}{\big|}_{\psi=0} = \varrho|_{\psi=0}\,.
\tag 2
$$
For terms undefined in this paper, cf., e.g., \cite{BF}, and \cite{S}.


\head 3. Theorem \endhead

\newpage

\proclaim{Theorem} Let $d \ge 3$ and $2 \le k \le d-1$ be integers, and
let $\lambda_0 \in (0, \infty)$. 
Suppose that for each $\lambda \in [0, \lambda _0]$, we have that
$K^{\lambda }$ is a convex
body in ${\Bbb{R}}^d$ 
with radial function $\varrho^\lambda (\omega)$, for
$\omega \in S^{d-1}$. Let $\varrho^0(\omega) \equiv 1$, and let
$\varrho^\lambda(\omega)$ be a $C^2$-function of $(\lambda,\omega) \in
[0,\lambda_0] \times S^{d-1}$. 
Assume that for each $\lambda \in [0, \lambda _0]$, 
for any linear
$k$-subspace $L_k \subset {\Bbb{R}}^d$, 
the function $y \mapsto W_1 (K^\lambda \cap
(L_k + y))$ has a maximum at $y=0$. Then
$$
\frac{\partial \varrho^\lambda}{\partial \lambda} (\omega){\big|}_{\lambda =0}
$$
is an even function of $\omega$.
If $2 \le l \le k-1$, and $\varrho^\lambda (\omega)$ is a
$C^3$-function of
$(\lambda,\omega) \in [0,\lambda_0] \times S^{d-1}$, and we
replace in the above hypothesis $W_1$ by $W_l$, then the same conclusion
holds.
\endproclaim

Clearly, we could have written, in the hypothesis of the theorem, that
$K^\lambda$ is a
star body, since, by the other assumptions, $K^\lambda$ is a convex body with
positive Gauss curvature for each $\lambda \in [0,\lambda_0]$ (after possibly
decreasing $\lambda_0$).

We observe that for the case $k=1$, and for the case $l=0$, we have the
theorems cited in the abstract, cf. \cite{H}, Theorem 1, and \cite{MM\'O},
Corollary 3.2. These assert
that, in this case, actually each $K^\lambda$ is 
centered, which
is of course a
stronger statement than the statement of the theorem of this paper. There is
still one
quermassintegral, namely $W_k$. However, this is, independently of its
argument, equal to the volume of the unit ball in ${\Bbb{R}}^k$, 
so, in this case the
hypotheses of our theorem do not imply anything.


\head 4. Proof \endhead

We begin with the following

\proclaim{Proposition} (\cite{MM\'O}, Theorem 3.8). Let $f: S^{d-1} \to 
{\Bbb{R}}$
be a $C^1$-function (or, more generally, a Lipschitz function). Further, let,
for each $\xi \in
S^{d-1}$ (or, more generally, for almost all $\xi \in S^{d-1}$), the equality
$$
\int \limits_{\xi^\bot \cap S^{d-1}} \frac{\partial f}{\partial \psi}
(\eta,0) d \eta = 0
$$
hold. Then $f$ is an even function. $\blacksquare $
\endproclaim


\demo{Proof of the Theorem}
As in the proof of Theorem 3.1 of \cite{MM\'O},
we
may suppose $k = d-1$. In fact, for any linear $(k+1)$-subspace $L_{k+1}$ of
${\Bbb{R}}^d$, we have
that $K^\lambda \cap L_{k+1}$ also satisfy the hypotheses of the theorem.
Furthermore, if for each $L_{k+1}$, the restriction of the function
$$
\frac{\partial \varrho^\lambda}{\partial \lambda} (\omega){\big|}_{\lambda = 0}
$$
to $L_{k+1}$ is even, then also this function itself is even. So, from now on,
let $k=d-1$.

Let $t_0 \in (0,\infty)$ be so small that the
closed ball about $0$, of radius $t_0$, is contained in each $K^\lambda$,
where $\lambda \in [0,\lambda_0]$ (possibly decreasing $\lambda_0$). From now
on, let

\newpage

$t \in (-t_0,t_0)$. This implies that $K \cap (\xi^\bot + t \xi)$ is a
($(d-1)$-dimensional) convex body in $\xi^\bot + t \xi$.

{\bf{1.}} First we 
treat
the case $l=1$. 

Let us fix a point $\xi \in
S^{d-1}$, that we consider as the north pole. Let $S^\lambda (\xi,t)$ denote
the ($(d-2)$-dimensional)
surface area of $K^\lambda \cap (\xi^\bot + t \xi)$, considered as a
($(d-1)$-dimensional) convex body in $\xi^\bot + t\xi$. We have
$$
S^\lambda (\xi,t) = \int \limits_{\xi^\bot \cap S^{d-1}} d S^\lambda (\xi,t)
= \int \limits_{\xi^\bot \cap S^{d-1}} r^\lambda(\xi,\eta,t)^{d-2}
\frac{1}{\langle \eta,n^\lambda(\eta,t)\rangle} d \eta\,,
\tag 3
$$
where $r^\lambda(\xi, \eta ,t)$
is the radial function of $K^\lambda \cap (\xi^\bot
+ t \xi)$, with respect to the ``origin'' $t \xi$, and $n^\lambda (\eta,t) \in
\xi^\bot \cap S^{d-1}$
is the outer normal unit vector of the surface element $d S^\lambda (\xi,t)$
at $\eta \in \xi^\bot \cap S^{d-1}$, taken in $\xi^\bot + t \xi$.

From now on, we
consider $\lambda \in [0,\lambda_0]$ as fixed, and drop the upper indices
$\lambda$. Also, to simplify the formulas, we omit those variables of our
functions, whose omission
does not lead to misunderstanding.

We determine
$$
\frac{\partial}{\partial t} S(\xi,t) |_{t=0}\,,
$$
that equals $0$ by the hypothesis of the theorem. We may
differentiate under
the integral sign. We have
$$
\cases
(\partial/\partial t) \left(r^{d-2} / 
\langle \eta, n (\eta)\rangle \right) =\\
(d-2) r^{d-3} (\partial r/\partial t) /
\langle \eta,n(\eta)\rangle
-r^{d-2} 
[(\partial /
\partial t) \langle \eta, n (\eta)\rangle ]/
\langle \eta,n(\eta)\rangle^2
\,,
\endcases
\tag 4
$$
and we have to evaluate this at $t=0$. 

Letting $t=0$, i.e., by \thetag{1},
$\psi =0$, we have by \thetag{2} $\partial t/\partial \psi = \varrho$, hence
$$
\frac{\partial}{\partial t} = \frac{\partial \psi}{\partial t} 
\frac{\partial}{\partial \psi} = \frac{1}{\varrho} \frac{\partial}
{\partial \psi}\,.
$$
Therefore, \thetag{4} equals
$$
(d-2) r^{d-3} \frac{1}{\varrho} \frac{\partial r}{\partial \psi} \frac{1}
{\langle \eta,n(\eta)\rangle} - r^{d-2}
\frac{1}{\langle \eta,n(\eta)\rangle^2} \frac{1}{\varrho}
\frac{\partial}{\partial \psi} \langle \eta, n(\eta)\rangle\,.
$$
Here the first term is, using $r = \varrho $ (cf. \thetag{1}),
$$
(d-2)\varrho ^{d-4}
\frac{\partial \varrho}{\partial \psi} \frac{1}{\langle
\eta, n(\eta)\rangle}\,,
\tag 5
$$
and the second term is
$$
-\varrho^{d-3} \frac{1}{\langle \eta,n(\eta)\rangle^2} \frac{\partial}
{\partial \psi} \langle n, n (\eta)\rangle\,.
\tag 6
$$

\newpage

Now it will be convenient to 
write
$\varrho = : 1 +
\varepsilon$, where $\varepsilon: S^{d-1} \to {\Bbb{R}}$ is a $C^2$-function, of
$C^2$-norm tending to $0$ for
$\lambda \to 0$. We calculate \thetag{5} and \thetag{6}, till terms of degree
$1$ in $\varepsilon$, but neglecting terms of degree at least $2$ in
$\varepsilon$.

Then \thetag{5}
becomes
$$
(d-2) (1+\varepsilon)^{d-4} \frac{\partial \varepsilon}{\partial \psi}
\frac{1}{\langle \eta, n(\eta)\rangle}\,.
$$
Here, because of the third factor, we may write $\varepsilon = 0$ in the
second and fourth factors, getting
$$
(d-2) \frac{\partial \varepsilon}{\partial \psi}\,.
\tag 7
$$
On the other hand, \thetag{6} contains $(\partial/\partial \psi) \langle \eta,
n(\eta)\rangle$ as a factor. We are going to show that this is an expression
of second 
order
in $\varepsilon$. We have
$$
n(\eta) = \frac{\left( 1,-\partial \varepsilon /\partial x_1,\dots,
-\partial \varepsilon /\partial x_{d-1}\right)}{\sqrt{1+( \partial
\varepsilon /\partial x_1)^2
+ \dots + \left( \partial \varepsilon /\partial x_{d-1} \right) ^2}}\,,
$$
where $x_1,\dots,x_{d-1}$ are the coordinates on $S^{d-1}$, in a neighbourhood
of $\eta$, given by the inverse of the exponential map at $\eta \in S^{d-1}$.
(The exponential map maps
vectors $u$, in a neighbourhood of the origin $\eta $ 
of the tangent plane of $S^{d-1}$ at
$\eta$, to the point $\omega \in S^{d-1}$ of the geodesic on $S^{d-1}$,
starting from $\eta$, in the
direction of $u$, with $\omega $ being
at a geodesic distance $\|u\|$ from $\eta$.) Therefore,
$$
\langle \eta, n(\eta) \rangle = \frac{1}{\sqrt{1+\left( \partial
\varepsilon /\partial x_1 \right)^2 + \dots +
\left(\partial \varepsilon /\partial x_{d-1} \right) ^2}}\,.
$$
Clearly, it is enough to show that, e.g.,
$$
\frac{\partial}{\partial x_1} \frac{1}{\sqrt{1+\left( \partial
\varepsilon /\partial x_1 \right)^2 + \dots + \left(\partial
\varepsilon /\partial x_{d-1} \right) ^2}}
\tag 8
$$
is of second degree of smallness in $\varepsilon$. However, \thetag{8} equals
$$
- \frac{(\partial \varepsilon /\partial x_1)  (\partial^2
\varepsilon /\partial x^2_1) + \dots + (\partial \varepsilon /\partial
x_{d-1})
\left( \partial^2\varepsilon /(\partial x_{d-1}\partial x_1) \right)
}{\left( 1+
\left(\partial \varepsilon /\partial x_1 \right)^ 2 + \dots +
\left(\partial \varepsilon /\partial x_{d-1} \right) ^2 \right) ^{3/2}}\,,
\tag 9
$$
and so our claim is shown.

Altogether, by \thetag{3} and \thetag{4}, 
and, on the one hand, by \thetag{5} and
\thetag{7}, on the other hand, by \thetag{6} and \thetag{9}, we have that
$((\partial/\partial t) S (\xi,t))|_{t=0}$ is, till terms of degree $1$ in
$\varepsilon$,
$$
(d-2) \int \limits_{\xi^\bot \cap S^{d-1}} \frac{\partial \varepsilon}
{\partial \psi} d \eta\,.
\tag 10
$$

\newpage

Since, for each $\xi \in S^{d-1}$,
\thetag{10} equals $0$, 
the Proposition implies that
$\varepsilon$ is even.
(Recall that, by hypothesis, $d \ge 3$.) Returning to
the original notations,
$$
\frac{\partial \varrho^\lambda}{\partial \lambda} (\omega) {\big|}_{\lambda=0}
$$
is an even function of $\omega$.

{\bf{2.}} 
Now we treat
the case $2 \le l \le d-2$.

Actually, we will allow $1 \le l \le d-1$. Of course, as stated after the
theorem, for $l=d-1$ the statement of the theorem does not hold. However, we
will need
this case for our formulas.

We have, for $1 \le l \le d-1$, that
$$
\cases
W^\lambda_l (\xi,t) : = W_l (K \cap (\xi^\bot + t\xi)) =\\
\left( 1/(d-1) \right) 
\int
\limits_{{\text{bd}}\,(K \cap(\xi^\bot + t\xi))} H_{l-1} (\xi,t) dS^\lambda
(\xi,t)\,.
\endcases
\tag 11
$$
Here $H_{l-1}(\xi,t)$ is
${\binom{d-2}{l-1}}^{-1}$ times the $(l-1)$'st elementary symmetric function of
the $d-2$ principal curvatures
$\kappa_1(\xi,t),\dots,\kappa_{d-2}(\xi,t)$ of 
bd\,$(K \cap (\xi^\bot + t\xi))$.
Cf., e.g., \cite{S}, p. 291.

We write $\kappa_i (\xi,t) = : 1 + \delta_i(\xi,t)$, where $\delta_i$ is of
first order
with respect to the $C^2$-norm of
$\varepsilon=\varrho - 1$.

Letting
$$
\cases
P: = {\binom{d-2}{l-1}}^{-1} \sum \limits_{1 \le i_1 < \dots < i_{l-1}\le
d-2} (1 + \delta_{i_1}) \dots (1 + \delta_{i_{l-1}})\\
-{\binom{d-2}{l-1}}^{-1} \sum \limits_{1 \le i_1 < \dots < i_{l-1}\le d-2}
(1 + \delta_{i_1} + \dots + \delta_{i_{l-1}})\\
={\binom{d-2}{l-1}}^{-1}
\sum \limits_{1 \le i_1 < \dots < i_{l-1}\le d-2} (1 + \delta_{i_1})
\dots (1 + \delta_{i_{l-1}})\\
- 1 - \left( (l-1)/(d-2) \right)
\sum \limits^{d-2}_{i=1}
\delta_i\,,
\endcases
\tag 12
$$
we have that $P$ is a linear combination with constant
coefficients, of the elementary symmetric functions of the $\delta_i$'s, of
degrees $2$ to $l-1$. Therefore, $\partial P/\partial t$ is a sum, whose
summands are
products of some $\partial \delta_i/\partial t$, and at least one further
$\delta_j$. Here $\partial \delta_i/\partial t$ is bounded by the
$C^3$-assumption, and the
$\delta_j$'s are 
of first order
with respect to the $C^2$-norm
of $\varepsilon$. Hence, when calculating the derivative of
\thetag{11}, with respect to $t$, at $t=0$, we can neglect
$\partial P/\partial t$. Hence, we
may replace in \thetag{11} $H_{l-1}$ by
$$
1 + \frac{l-1}{d-2} \sum \limits^{d-2}_{i=0} \delta_i\,,
$$
and this replacement will not affect the calculation of the
derivative of \thetag{11}, with respect to $t$, at $t=0$.

\newpage

We turn to the calculation of the derivative of \thetag{11}, with
respect to $t$, at
$t=0$, which has to be $0$. As mentioned above,
this equals
$$
\cases
\left(1/(d-1) \right) (\partial/\partial t) \int \limits_{{\text{bd}}\,
(K \cap (\xi^\bot + t\xi))} dS^\lambda (\xi,t) +\\
(l-1)/\left( (d-1)(d-2) \right)
(\partial/\partial t)
\int \limits_{{\text{bd}}\,(K \cap (\xi^\bot + t \xi))}
\left(\sum \limits^{d-2}_{i=0} \delta_i\right) dS^\lambda (\xi,t)\,.
\endcases
\tag 13
$$
Here the first summand is, by {\bf{1}},
$$
\frac{d-2}{d-1} \int \limits_{\xi^\bot \cap S^{d-1}} \frac{\partial
\varepsilon}{\partial \psi} d \eta\,.
$$
We are going to determine the second summand. For this, put $l=d-1$. Then, as
already mentioned, \thetag{11} is constant, hence \thetag{13} equals $0$. From
this we have
$$
\frac{\partial}{\partial t} \int \limits_{{\text{bd}}\,
(K \cap (\xi^\bot + t\xi))} \left(\sum \limits^{d-2}_{i=0}
\delta_i\right) dS^\lambda (\xi,t) = -(d-2)
\int \limits_{\xi^\bot \cap S^{d-1}} \frac{\partial \varepsilon}
{\partial \psi} d\eta\,.
$$
Hence, for all $l = 1, \dots, d-1$, we have that \thetag{13} further equals
$$
\frac{d-1-l}{d-l} \int \limits_{\xi^\bot \cap S^{d-1}} \frac{\partial
\varepsilon}{\partial \psi} d \eta\,,
\tag 14
$$
which
equals $0$. By the hypothesis of the theorem, we have 
$l \le d-2$, hence
$$
\int \limits_{\xi^\bot \cap S^{d-1}} \frac{\partial \varepsilon}
{\partial \psi} d \eta = 0\,,
$$
for each $\xi \in S^{d-1}$. As in {\bf{1}}, this implies that
$$
\frac{\partial \varrho^\lambda}{\partial \lambda} (\omega) |_{\lambda=0}
$$
is an even function of $\omega$. 
$\blacksquare $
\enddemo


\head 5. Remark \endhead

\definition{Remark} Let $\lambda_0 \in (0,\infty)$. Let $K^0 \subset
{\Bbb{R}}^d$ be a centered convex body, further
suppose that for each $\lambda \in (0,\lambda_0]$, we have that
$K^{\lambda }$ is a convex
body in ${\Bbb{R}}^d$, with radial functions
$\varrho^\lambda$, for $\lambda \in
\{ 0 \} \cup (0, \lambda _0]=[0,\lambda_0]$. 
Moreover, let $\varrho^\lambda (\omega)$ be a $C^2$-function
of $(\lambda,\omega)
\in [0,\lambda_0] \times S^{d-1}$. We may ask whether some analogue of our

\newpage

theorem holds. That is, suppose that for each $\lambda \in [0,\lambda_0]$,
and each linear $(d-1)$-subspace
$L_{d-1} \subset {\Bbb{R}}^d$, the function $y \mapsto
W_l(K^\lambda \cap (L_{d-1} + y))$
has a maximum at $y=0$. Then we may pose the question: is
$$
\frac{\partial \varrho^\lambda}{\partial \lambda} (\omega) |_{\lambda=0}
$$
an even function of $\omega$? However, we will show that this question, even
in the simplest unsolved case, i.e., for $d=3$, and for $W_1$, is untreatable
by our present
methods.
\enddefinition

For $d-1=2$ we can use, for the calculation of the perimeter of
$K^\lambda \cap (\xi^\bot + t \xi)$, the simpler formula
$ds^2 = dr^2 + r^2 d\eta^2$. Then the equality
$$
\left(\frac{\partial}{\partial t} W_1 \left[(K \cap (\xi^\bot + t \xi))-t
\xi\right]\right){\big|}_{t=0} = 0
$$
can be rewritten as
$$
\int \limits_{S^1} \frac{1}{\sqrt{\varrho^2+
(\partial \varrho /\partial \eta )^2}} \left(\frac{\partial \varrho}
{\partial \psi} + \frac{1}{\varrho} \frac{\partial \varrho}{\partial \eta}
\frac{\partial^2 \varrho}{\partial \eta \partial \psi}\right) d\eta = 0\,.
\tag 15
$$
Let us write $\varrho^\lambda = \varrho^0 + \varepsilon$. We retain in 
\thetag{15} the
terms at most linear in $\varepsilon$, 
and investigate this situation.
Clearly, the terms of degree $0$ in $\varepsilon $ 
together give the
integral, on $S^1$, of an odd function, i.e., $0$. Now we investigate the
terms of degree $1$ in $\varepsilon $, in the expression
under the integral sign in \thetag{15}. These are the following: 
$$
\cases
\left[-
[(\varrho ^0)^2 +
(\partial{\varrho^0}/\partial \eta)^2]^{-1/2} 
(\varrho^0)^{-2}
(\partial{\varrho^0}/\partial \eta)
\partial^2{\varrho^0}/(\partial\eta \partial\psi )\right.\\
- 
{\varrho^0 ((\varrho^0)^2 +
(\partial{\varrho^0}/\partial\eta)^2)^{-3/2}}
\left.\left(
\partial{\varrho^0}/\partial\psi +
(\varrho^0)^{-1}
(\partial \varrho ^0/ \partial \eta )
\partial ^2 {\varrho^0}/(\partial\eta
\partial\psi)\right)\right] \varepsilon  \\
+
[(\varrho^0)^2 + (\partial \varrho^0 /\partial\eta)^2]^{-1/2}
(\partial\varepsilon/\partial\psi) +
\left[
[(\varrho^0)^2 +
(\partial{\varrho^0}/\partial\eta)^2] ^{-1/2}
(\varrho^0) ^{-1}
\partial^2 \varrho^0/\right.\\
\left.
(\partial\eta \partial\psi) -
(\partial{\varrho^0}/\partial\eta )[(\varrho^0)^2 +
(\partial {\varrho^0}/\partial\eta)^2]^{-3/2}
\left(
\partial \varrho^0/\partial \psi +
(\varrho^0)^{-1} 
(\partial{\varrho^0}/\partial \eta)
\partial^2\varrho^0/
\right. \right.
\\
\left. \left.
(\partial\eta \partial\psi)\right)\right]
\partial \varepsilon/\partial \eta
+ 
[(\varrho^0)^2 + (\partial \varrho^0/\partial\eta)^2]^{-1/2}
(\varrho^0)^{-1} 
(\partial \varrho^0/\partial\eta)
\partial^2\varepsilon/(\partial\eta \partial\psi)\\
=: A \, \varepsilon +B \, \partial \varepsilon / \partial
\psi +C \, \partial \varepsilon /\partial \eta  +D \, \partial ^2
\varepsilon /(\partial \eta \partial \psi )\,. 
\endcases
\tag 16
$$
Now, let us suppose that $\varrho ^{\lambda }(\omega )$ is a $C^3$ function of 
$(\lambda,\omega) \in [0,\lambda_0] \times S^2$.
Then,
retaining in 
\thetag{15} the
terms at most linear in $\varepsilon $,
\thetag{15} becomes, by integration by parts,
$$
\int _{S^1} \left[ \left( A-\frac{\partial C}{\partial \eta } \right)
\varepsilon + \left( B-\frac{\partial D}{\partial \eta } \right)
\frac{\partial \varepsilon }{\partial \psi } \right] d \eta =0\,.
\tag 17
$$
(We do not give the coefficients in this formula more explicitly.) 
Of course, the left hand
side of \thetag{17} is a
continuous linear operator in $\varepsilon $, for the $C^1$-topology. But its
solution (e.g., that the solutions among the
$C^1$-functions would be just the even $C^1$-functions)
seems to be untreatable by our methods.


%
%
%
%
%
%

\enddocument